\documentclass[graybox]{svmult}

\usepackage{type1cm}        
\usepackage{makeidx}         
\usepackage{graphicx}        
\usepackage{multicol}        
\usepackage[bottom]{footmisc}
\usepackage{placeins}
\usepackage{subcaption}
\usepackage{newtxtext}       %
\usepackage[varvw]{newtxmath}       

\usepackage{multicol}
\usepackage{multirow}

\def\bA{{\bf A}}
\def\bR{{\bf R}}

\def\bD{{\bf D}}

\def\dsp{\displaystyle}

\def\O{\Omega}

\def\G{\Gamma}

\def\to{\rightarrow}

\def\R{\mathbb{R}}

\def\H10{H^1_{ \partial \O \setminus \partial \O_S}}

\def\dsp{\displaystyle}

\def\nodes{{\cal V}}

\makeindex             

\begin{document}
	
	\title*{A Trefftz-like coarse space for the two-level Schwarz method on perforated domains}
	\author{Miranda Boutilier, Konstantin Brenner, and Victorita Dolean}
	\institute{Miranda Boutilier \at Universit{\'e} C{\^o}te d'Azur, LJAD, \email{miranda.boutilier@univ-cotedazur.fr}
		\and Konstantin Brenner \at Universit{\'e} C{\^o}te d'Azur, LJAD, CNRS, INRIA, \email{konstantin.brenner@univ-cotedazur.fr}
		\and Victorita Dolean \at University of Strathclyde, Dept. of Maths and Stats and Universit{\'e} C{\^o}te d'Azur, LJAD, CNRS \email{work@victoritadolean.com}}
	%
	%
	\maketitle
	
	\abstract*{			We consider a new coarse space for the ASM and RAS preconditioners to solve elliptic partial differential equations on perforated domains, where the numerous polygonal perforations represent structures such as walls and buildings in urban data. 
		With the eventual goal of modelling urban floods by means of the nonlinear Diffusive Wave equation, this contribution focuses on the solution of linear problems on perforated domains.
		Our coarse space uses a polygonal subdomain partitioning and is spanned by Trefftz-like basis functions that are piecewise linear on the boundary of a subdomain and harmonic inside it. It is based on nodal degrees of freedom that account for 
		the intersection between the perforations and the subdomain boundaries. As a reference, we compare this coarse space to the well-studied Nicolaides coarse space with the same subdomain partitioning. It is known that the Nicolaides space is unable to prevent stagnation in convergence when the subdomains are not connected; we work around this issue by separating each subdomain by disconnected component. Scalability and robustness are tested for multiple data sets based on realistic urban topography. Numerical results
		show that the new coarse space is very robust and accelerates the number of
		Krylov iterations when compared to Nicolaides, independent of the complexity of the data.}
	

	\section{Introduction and Model Problem}
	\label{sec:1}
	
	
	Numerical modeling of overland flows plays an increasingly important role in predicting, anticipating and controlling floods, helping to size and position protective systems including dams, dikes or rainwater drainage networks. One of the challenges of the numerical modeling of urban floods is that the small structural features (buildings, walls, etc.) may significantly affect the flow.  
	Luckily,  modern terrain survey techniques including photogrammetry and Laser Imaging, Detection, and Ranging (LIDAR) allow to acquire high-resolution topographic data for urban areas as well as for natural (highly vegetated) media.
	For example, the data set used in this article has been provided by Métropole Nice Côte d'Azur (MNCA) and allows for the infra-metric description of the urban geometries \cite{buildings}.  
	
	
	From the hydraulic perspective,  these structural features can be assumed to be essentially impervious,  and therefore represented as perforations (holes) in the model domain.  
	Our long term modelling strategy is based on the 
	Diffusive Wave equation \cite{diffwave}. However,
	understanding linear problems posed on perforated domains is a crucial preliminary step and the object of this contribution.  
	
	Let $D$ be an open simply connected polygonal domain in $\mathbb{R}^2$,  we denote by $\left(\Omega_{S,k}\right)_k$ a finite family of perforations in $D$ such that each $\Omega_{S,k}$ is an open connected polygonal subdomain of $D$. The perforations are mutually disjoint, that is $\overline{\Omega_{S,k}}\cap \overline{\Omega_{S,l}} = \emptyset$ for any $k\neq l$.  We denote $\Omega_S = \bigcup_k \Omega_{S,k}$ and $\Omega = D\setminus \overline{\Omega_S}$, assuming that the family $\left(\Omega_{S,k}\right)_k$ is such that $\Omega$ is connected.  Note that the latter assumption implies that $\Omega_{S,k}$ are simply connected.
	
	Let $f\in L^2(\Omega)$,  in this article we are interested in the boundary value problem
	\begin{equation}\label{model_pde}
		\left\{
		\begin{array}{rll}
			- \Delta u &=& f \qquad \mbox{in} \qquad \Omega, \\
			\dsp - \frac{\partial u}{\partial \bf{n}} &=& 0 \qquad \mbox{on} \qquad \partial  \Omega  \cap \partial \Omega_S,\\
			u &=& 0 \qquad \mbox{on} \qquad  \partial \Omega \setminus \partial \Omega_S.\\
		\end{array}
		\right.
	\end{equation}
	
	Depending on the geometrical complexity of the computational domain, the numerical resolution of \eqref{model_pde} may become challenging. A typical data set that we are interested in, illustrated by Figure \ref{fig:a}, may contain numerous perforations that are described on different scales. In this regard, our strategy relies on the use of a Krylov solver combined with domain decomposition (DD) methods. Generally, to achieve scalability with respect to the number of subdomains in overlapping Schwarz methods, coarse spaces/components are needed. Including a coarse space in a Schwarz preconditioner results in what is referred to as a two-level Schwarz preconditioner.
	
	The model problem can be thought of as the extreme limit case of the elliptic model containing highly contrasting coefficients. 
	Two-level domain decomposition methods have been extensively studied for such heterogeneous problems.  
	There are many classical results for coarse spaces that are contructed so as to resolve the jumps of the coefficients; see \cite{analysis, multileveldiscon, balancing}   for further details.
	Approaches to obtain a robust coarse space without careful partitioning of the subdomains include spectral coarse spaces such as those given in \cite{eigenproborig,dtn, geneo}.  Additionally, the family of GDSW
	(Generalized Dryja, Smith, Widlund) 
	methods \cite{gdsw}  employ energy-minimizing coarse spaces and can be used to solve heterogeneous problems on less regular domains. These spaces are discrete in nature and involve both edge and nodal basis functions.

Alternatively, robust coarse spaces can be constructed using the ideas from
multi-scale finite elements methods (MsFEM) \cite{multiscalehighaspectratio, multiscalepdes}. 
The combination of spectral and MsFEM methods can be found in \cite{shemorig}. Outside of the DD framework, 
specifically on domains with small and numerous perforations, the authors of \cite{legoll, taralova} also introduced an enriched MsFEM-like method. 

Here, we present an efficient and novel coarse space in the overlapping Schwarz framework inspired by the Boundary Element based Finite Element (BEM-FEM) method \cite{bemfem}. In contrast with the classical BEM-FEM approach, the local multiscale basis functions are computed numerically such as in MsFEM methods.  
This approach is motivated by our interest in nonlinear time dependent models for which the analytical expression of the fundamental solutions may not be easily available.


\section{Discretization and Preliminary Notations}\label{sec:asm}
We introduce a coarse discretization of $\Omega$ which involves a family of polygonal cells $\left( \Omega_j \right)_{j= 1,\ldots, N}$, the so-called coarse skeleton $\Gamma$, and the set of coarse grid nodes that will be referred to by $\nodes$.

The construction is as follows. Consider a finite nonoverlapping polygonal partitioning of $D$ denoted by $\left( D_j \right)_{j= 1,\ldots, N}$ and an induced nonoverlapping partitioning of $\Omega$ denoted by $\left( \Omega_j \right)_{j= 1,\ldots, N}$  such that $\Omega_j=D_j \cap \Omega$. We will refer to  $\left( \Omega_j \right)_{j= 1,\ldots, N}$  as the coarse mesh over $\Omega$.
Additionally,  we denote by $\Gamma$ the skeleton of the coarse mesh, that is  $\Gamma=  \bigcup_{j \in \{1, \ldots, N\} } \partial \Omega_j \setminus \partial\Omega_S$.

\begin{figure}
	\sidecaption[t]
	\centering
	\includegraphics[scale=.74]{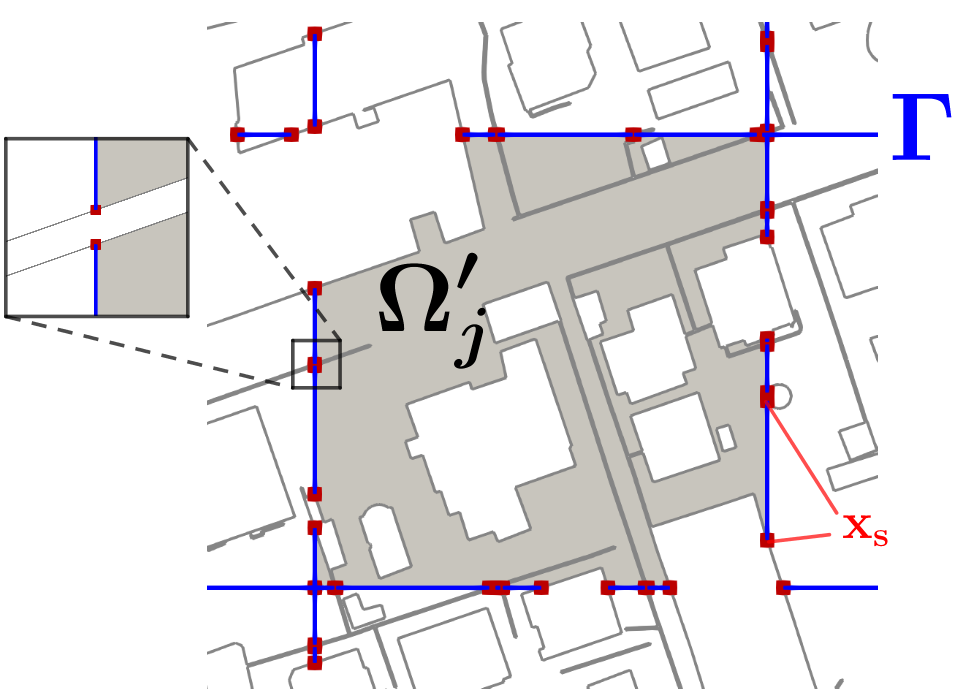}
	\newline
	\caption{ Coarse grid cell $\Omega_j$, nonoverlapping skeleton $\Gamma$ (blue lines), and coarse grid nodes $\mathbf{x}_s \in \nodes$ (red dots). Coarse grid nodes are located at $\overline{\Gamma} \cap \partial \Omega_S.$}
	\label{fig:coarsedofs}
\end{figure}

Let $\text{vert}(\Omega_j)$ denote the set of vertices of the polygonal domain $\Omega_j$. The set of coarse grid nodes is given as 
$\nodes= \bigcup_{j \in \{1, \ldots, N\} }   \text{vert}(\Omega_j) \cap \overline{\Gamma}$.  The total number of coarse grid nodes is denoted by $N_\nodes$.
We refer to Figure  \ref{fig:coarsedofs} for the illustration of the coarse mesh entities.

\begin{figure}
	\centering
	\begin{subfigure}{0.49\textwidth}
		\centering
		\includegraphics[width=.8\linewidth]{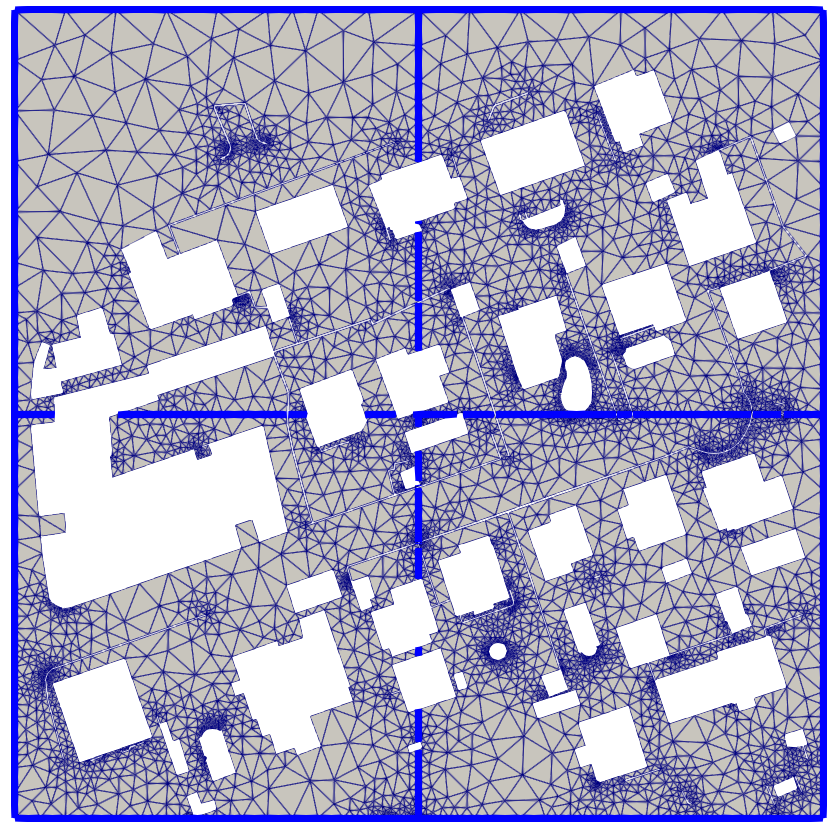}
		\caption{2$\times$2 subdomains}
	\end{subfigure}
	\begin{subfigure}{0.49\textwidth}
		\centering
		\includegraphics[width=.8\linewidth]{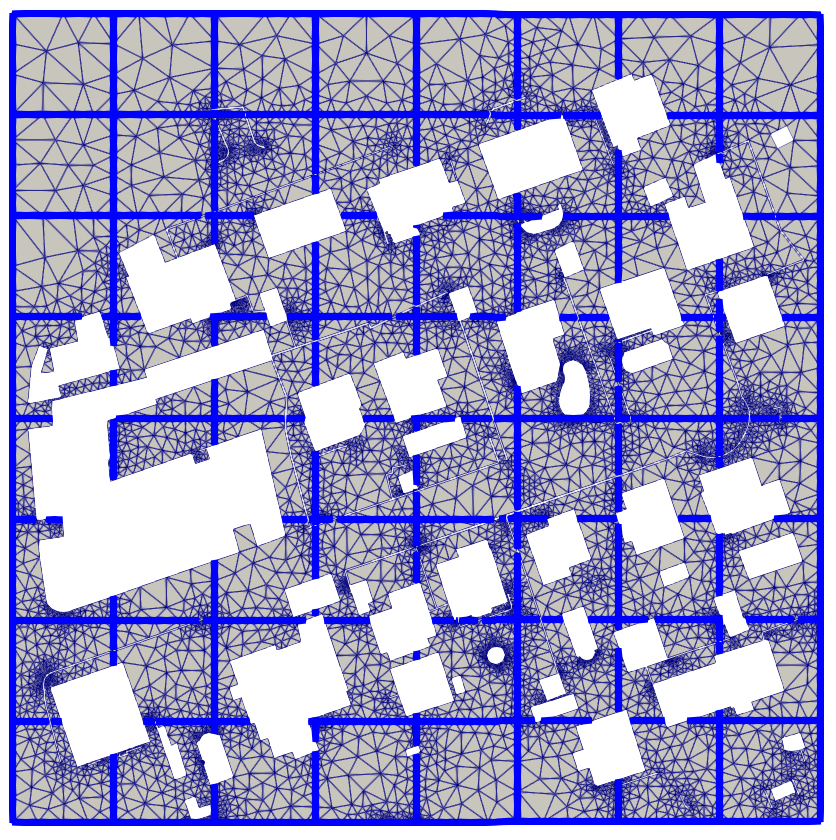}
		\caption{8$\times$8 subdomains}
	\end{subfigure}
	\caption{Conforming triangulation for the same domain $\Omega$ with different numbers of coarse cells $N$. The coarse skeleton $\Gamma$ is shown by the blue lines.}
	\label{fig:triang}
\end{figure}

We discretize the model problem \eqref{model_pde} with piecewise linear continuous finite elements on a triangular mesh of $\Omega$. This mesh is conforming to the coarse polygonal $\left( \Omega_j \right)_{j= 1,\ldots, N}$; an example of the triangulation for various numbers of coarse cells $N$ is given in Figure \ref{fig:triang}.
The finite element discretization of \eqref{model_pde} results in the linear system
$\bA \mathbf{u}=\mathbf{f}.$

Let $\left( \Omega_j'\right)_{j= 1,\ldots, N}$ denote the set of overlapping subdomains of $\Omega$.  In practice, each $\Omega_j'$ is constructed by propagating $\Omega_j$ by a few layers of triangles. 
Consider classical boolean restriction matrices $\bR_j$ and corresponding extension matrices  $\bR_j^T$ associated to the family of overlapping subdomains $\left( \Omega_j'\right)_{j= 1,\ldots, N}$.
With a coarse restriction matrix $\bR_0$ that will be specified later, the two-level discrete Additive Schwarz (ASM) preconditioner is given by 
\begin{equation}\label{eq:RAS2}
	M_{ASM,2} ^{-1}= \bR_0^T(\bR_0\bA \bR_0^T)^{-1}\bR_0+ \sum_{j=1}^{N} \bR_j^T\b (\bR_j\bA\bR_j^T)^{-1}\bR_j.
\end{equation}

\section{Description of the Trefftz-like coarse space}\label{sec:trefftz}

Here we introduce the  Trefftz-like coarse space spanned by the functions that are piecewise linear on the skeleton $\G$ and discrete harmonic inside the nonoverlapping subdomains
$\Omega_j$.   
For any node  $\mathbf{x}_s \in  \nodes$,  we introduce the function $g_s: \Gamma \rightarrow \R$,  which is
continuous on $\G$ and linear on each edge of $\Gamma$.  It is clear that $g_s$ is fully defined by its values at the nodes $\mathbf{x}_i \in \nodes$, for which we set
\[ g_s(\mathbf{x}_i)=\begin{cases} 
	1, & s=i,   \\
	0, & s \neq i. \\
\end{cases}
\]
To illustrate the construction of the nodal basis of the coarse space,  we consider the following set of boundary value problems. For all $\Omega_j$ and for all $s=1, \ldots, N_\nodes$, find  $\phi_{s}^j \in H^1(\Omega_j)$ such that $\phi_s^j$ is the weak solution to the following problem
\begin{equation}\label{eq:laplacedir}
	\left\{
	\begin{array}{rllll}
		- \Delta \phi_{s}^j &=& 0  &  \mbox{in} &  \Omega_j, \\
		\dsp -	\frac{\partial \phi_{s}^j}{\partial \bf{n}}&=&0 &\mbox{on} & \partial \Omega_j \cap \partial \Omega_S,\\
		\phi_{s}^j &=& g_s &\mbox{on} &  \partial \Omega_j\setminus \partial \Omega_S.\\
	\end{array}
	\right.
\end{equation}
The finite element discretization of \eqref{eq:laplacedir} results in the system of the form 
$
\bA_j' \boldsymbol{\phi}_{s}^j=\mathbf{b}_s^j,
$
where $\bA_j' $ is the local stiffness matrix and $\mathbf{b}_s^j$ accounts for the Dirichlet boundary data in \eqref{eq:laplacedir}.
Let $\overline{\bR}_j$ denote the restriction matrices corresponding to $\overline{\Omega_j}$,  and let $\boldsymbol{\phi}_s$ be a vector such that
$
\overline{\bR}_j \boldsymbol{\phi}_s=  \boldsymbol{\phi}_{s}^j 
$
for all $j = 1,\ldots, N$.
The coarse space is then defined as the span of the basis functions $\phi_s,  s=1, \ldots, N_\nodes$, while the $k$th row of $\mathbf{R}_0$ is given by $\boldsymbol{\phi}_k^T$ for $k=1, \ldots, N_\nodes$.


\section{Numerical Results}

We present below the numerical experiments concerning the performance of the conjugate gradient (CG) method using  two-level preconditioner \eqref{eq:RAS2} and the Trefftz-like coarse space introduced in Section \ref{sec:trefftz}. For the sake of comparison we also report the numerical results obtained using a more standard Nicolaides coarse space which is going to be detailed later.

\begin{figure}
	\centering
	\begin{subfigure}{0.49\textwidth}
		\centering
		\includegraphics[width=.8\linewidth]{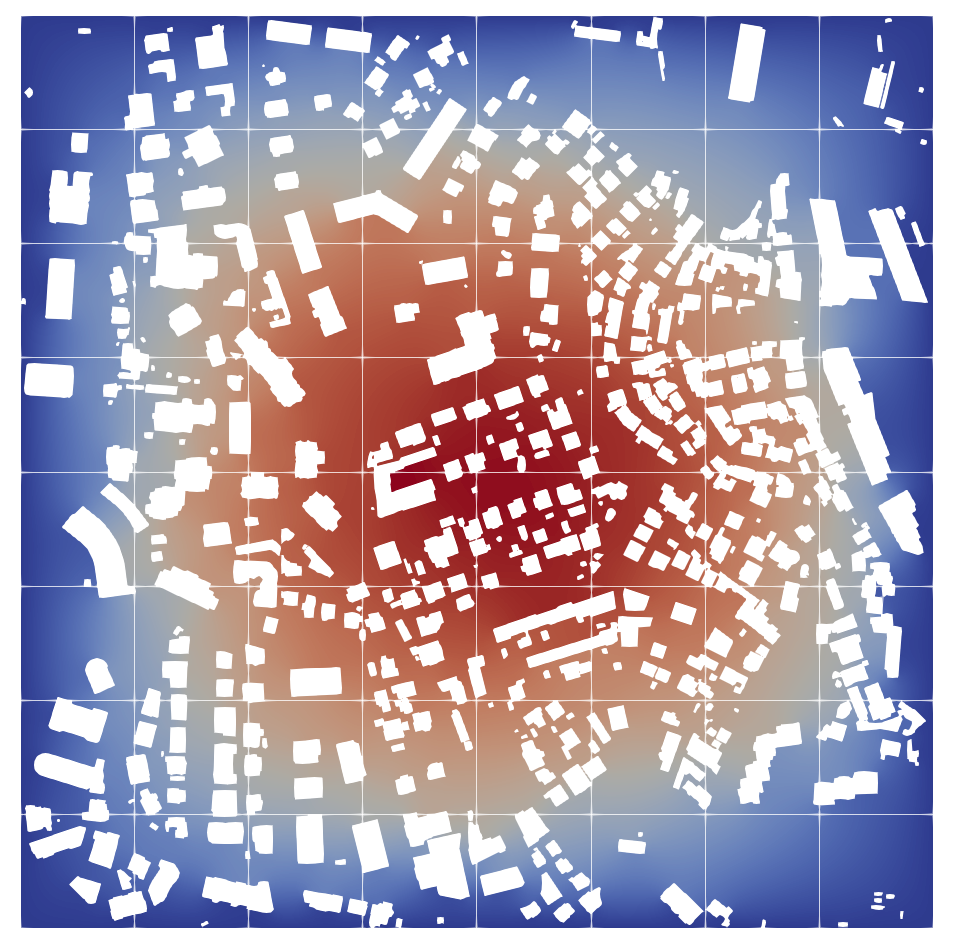}
		\caption{Without walls}
	\end{subfigure}
	\begin{subfigure}{0.49\textwidth}
		\centering
		\includegraphics[width=.8\linewidth]{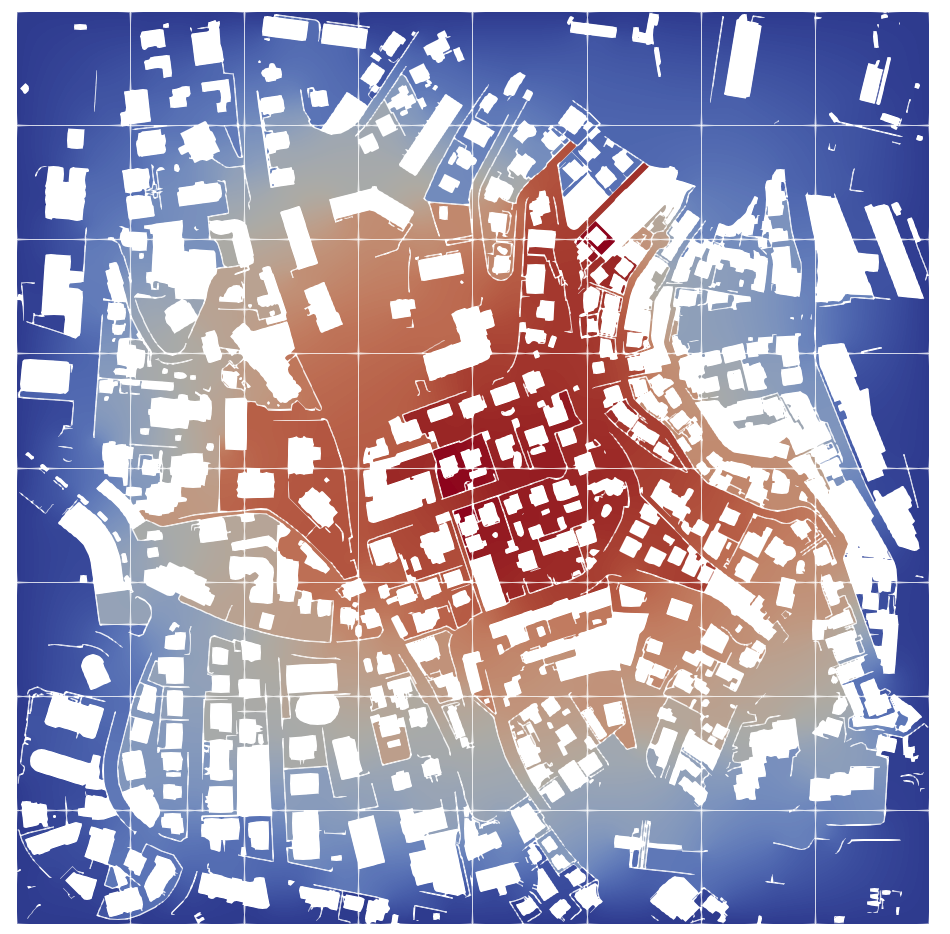}
		\caption{With walls}
	\end{subfigure}
	\caption{
		Approximate solution over a computational domain divided in $N=8 \times 8$ nonoverlapping subdomains.
}
\label{fig:a}
\end{figure}

The data sets used in this experiment have been kindly provided by Métropole
Nice Côte d'Azur  and reflect the structural topography of the city of Nice.  Although this type of data is available for the whole city \cite{buildings},  we focus here on a relatively small special frame (see Figure \ref{fig:a}).   In this numerical experiment  we consider two kinds of structural elements - buildings (and assimilated small elevated structures) and walls.  We note that the perforations resulting from the data sets we use (especially the wall data) can span across multiple coarse cells,  which is a challenging situation for traditional coarse spaces.


%

In this numerical experiment we consider the problem \eqref{model_pde} with the-right-hand side given by $f = 1$.  Figure \ref{fig:a} reports the finite element solution obtained for the data excluding and including walls. 
The figure also reflects the nonoverlapping partitioning into $N = 8\times 8$ subdomains.

\begin{figure}
\centering
\begin{subfigure}{0.495\textwidth}
\centering
\includegraphics[width=\linewidth]{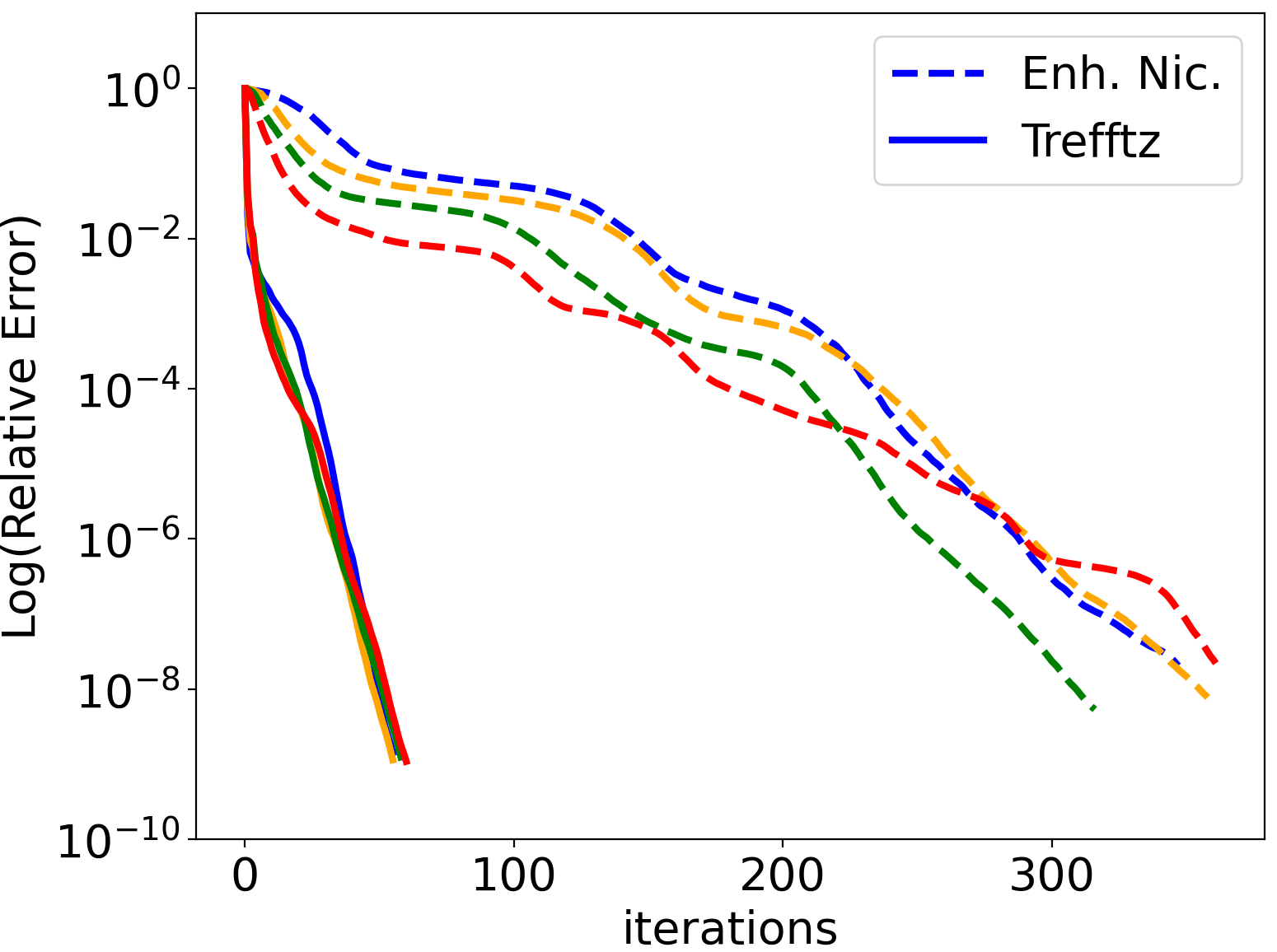}
\caption{Minimal geometric overlap}
\end{subfigure}
\begin{subfigure}{0.495\textwidth}
\centering
\includegraphics[width=\linewidth]{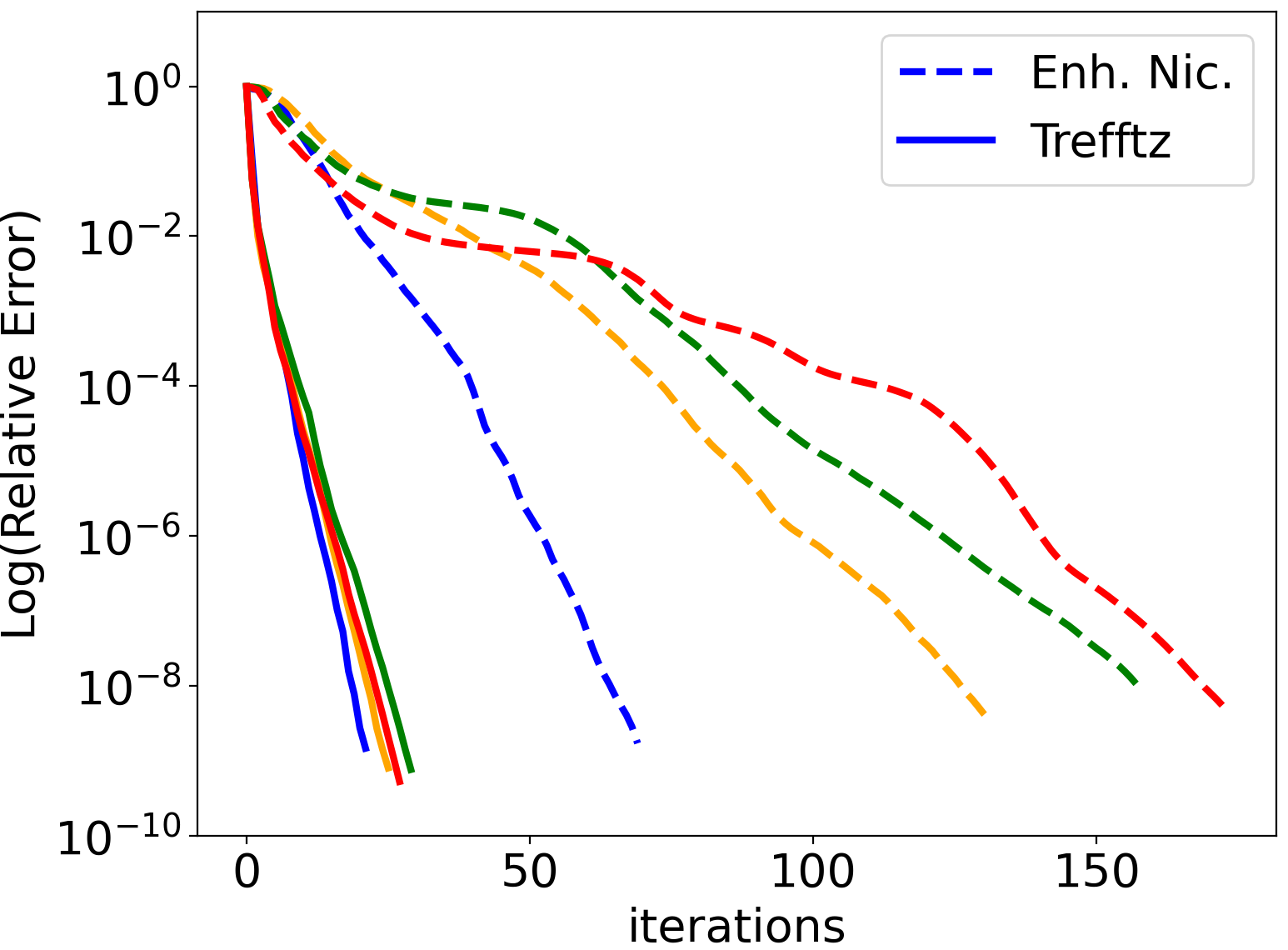}
\caption{Overlap $\frac{1}{20}L$}
\end{subfigure}
\caption{Convergence curves for the Trefftz-like (solid lines) and Nicolaides (dashed lines) coarse spaces for the data set involving both buildings and walls and two overlap sizes. Colors correspond to the number of subdomains as follows: $N=16$ (blue), $N=64$ (orange), $N=256$ (green), $N=1024$ (red).}
\label{fig:convcurves}
\end{figure}

Figure \ref{fig:convcurves} and Table \ref{tab:table} report the performance of the two-level preconditioner used in the PCG method, for varying number of subdomains $N$ and two relative overlap sizes.  As the computational domain $\Omega$ remains fixed independently of $N$,  the results of this experiment could be interpreted in terms of a strong scalability.  However, we wish to stress that the fine-scale triangulation is obtained based on the nonoverlapping partitioning $(\Omega_j)_{j = 1,\ldots, N}$.  Consquentially, the linear system $\bA \mathbf{u}=\mathbf{f}$ changes from one coarse partitioning to another.   Nevertheless we ensure that the dimension of  the system is roughly constant throughout the experiment.  Depending on the chosen $N$, the linear system involves about $60k$ (buildings alone)  and $180k$ (buildings and walls) nodal unknowns. 


For the sake of comparison, we also provide numerical results for the well-known Nicolaides coarse space \cite{nic}, made of flat-top partition of unity functions associated with the overlapping partitioning.  
As the scalability provided by the Nicolaides space relies on the Poincaré inequality over the subdomains, 
we further partition $(\Omega_j')_j$ into a family of connected regions for this space. 
In other words, let $m_j$ denote the number of disconnected components for each overlapping subdomain $\Omega_j'$
and let $\Omega_{j,l}', l=1, \ldots, m_j$ denote the corresponding disconnected component. Then our new overlapping partioning contains $m=\sum_{j=1}^{N} m_j$ total subdomains and is given by
$$ (\hat{\Omega}_{k})_{k \in \{1, \ldots, m\}} = \bigl( (\Omega_{j,l}')_{l \in \{1, \ldots, m_j\}} \bigr)_{j \in \{1, \ldots, N\}}. $$
Then, the Nicolaides coarse space is as follows. The $k$th row of $\mathbf{R}_0$ and therefore the $k$th column of $\bR_0^T$ is given by $(\mathbf{R_0^T})_k=  \widehat{\bR}_k^T \widehat{\bD}_k \widehat{\bR}_k\mathbf{1}$
for $k=1, \ldots, m$, where $\widehat{\bR}_k$ and $\widehat{\bD}_k$ are the restriction and partition of unity matrices corresponding to $\hat{\Omega}_{k}$ and $\mathbf{1}$ is a vector full of ones. The partition of unity matrices are constructed such that 
$\mathbf{I}= \sum_{k=1}^{m} \widehat{\bR}_k^T \widehat{\bD}_k \widehat{\bR}_k.$

Figure \ref{fig:convcurves} reports 
convergence histories of the preconditioned  CG method using the Nicolaides and Trefftz-like coarse spaces for the data set including both walls and buildings.
Table \ref{tab:table} summarizes the numerical performance for data sets including or excluding walls.  In particular, for both preconditioners,  it reports 
the dimensions of the coarse spaces, as well as the number of CG iterations required to achieve a relative $l^2$ error of $10^{-8}$.

The performance of the Trefftz-like coarse space appears to be very  robust with respect to both $N$ and the complexity of the computational domain.  The improvement with respect to the alternative Nicolaides approach is quite striking, especially in the case of the minimal geometric overlap.
As expected,  increased overlap in the first level of the Schwarz preconditioner provides additional acceleration in terms of iteration count.  However,  for the Trefftz-like space,  the results with minimal geometric overlap appear to already be quite reasonable. 

The dimensions and the relative dimensions of the two coarse spaces are reported in Table \ref{tab:table}.
Relative dimension refers to the would-be dimension of the coarse space in the case of a homogeneous domain with $\Omega_S = \emptyset$,  that is, the relative dimensions are computed as $\frac{\text{dim}(R_0)}{(\sqrt{N}+1)^2}$ for the Trefftz-like space and as $\frac{\text{dim}(R_0)}{N}$ for the Nicolaides space.  
We observe that the Trefftz-like coarse space requires a much larger number of degrees of freedom,
which naturally leads to a large coarse system to solve.   We note that the contrast between the dimensions of two spaces reduces as $N$ grows.  In general,  the dimension of the Trefftz-like coarse space seems reasonable given the geometrical complexity of the computational domain.

\begin{table}
\centering
\caption{ PCG iterations, condition number, dimension, and relative dimension for the Trefftz-like and Nicolaides coarse spaces.  Results are shown for minimal geometric overlap and $\frac{1}{20}H$, where $H=\max_{j} \text{diam}(\Omega_j)$. As the dimension of the Nicolaides space will change with respect to the overlap, its dimension is given as the average dimension over the two overlap values.
}
\setlength{\tabcolsep}{4.pt}
\begin{tabular}	{|c r |cc|cc|c|cc|cc|c|} 
\hline
& & \multicolumn{5}{|c}{Nicolaides}  &  \multicolumn{5}{|c|}{Trefftz}
\\ 
\hline
& &	\multicolumn{2}{|c|}{it.} & \multicolumn{2}{|c|}{cond.}  &  \multicolumn{1}{|c|}{dim. (rel)}	  &	\multicolumn{2}{|c|}{it.} & \multicolumn{2}{|c|}{cond.}  &  \multicolumn{1}{|c|}{dim. (rel)}	      \\
N	&&min. & $\frac{H}{20}$  & min. &  $\frac{H}{20}$  &&min &  $\frac{H}{20}$  & min. &  $\frac{H}{20}$ & \\
\hline
16& no walls&149& 51   &     581    &    82        & 21 (1.3)  &52&  28&   59    &       11                 &170 (6.8)    \\
& walls &348&     70     &     6826    &   133       &96 (6.0)  &56&  22 &   136    &        7               &400  (16.0)   \\
\hline
64& no walls&164&  78 &      567   &   119         &85 (1.3)  &50&  28  &   50  &             12         &433  (5.3)   \\
& walls &359& 132       &       5902    &   297    &256   (4.0) &56&26   &  57   &          9           &880  (10.9)   \\
\hline
256& no walls&136& 81&      273 &      89        &312  (1.2)  &56&   27&   54  &      10               &1010  (3.5)   \\
& walls &317&   159   &    4575     &          12    &719   (2.8)  &59 & 30 & 60    &               13      &1912  (6.6)  \\
\hline
1024& no walls&120&83&     341   & 149                 &1204 (1.2)&56&   28 &  76 &        13           &2500  (2.3)   \\
& walls &362&  174  &     3895   &    1310        &2044   (2.0)  &61& 28 &  97   &      13              &4253   (3.9)\\
\hline
\end{tabular}
\label{tab:table}
\end{table}

\section{Conclusions}
In this work we presented a novel Trefftz-like coarse space for the two-level ASM preconditioner,  specifically designed for problems resulting from  elliptic PDEs in perforated domains.
This coarse space is robust with respect to data complexity and number of subdomains on a fixed total domain size,  and provides significant acceleration in terms of Krylov iteration counts when compared to a more standard Nicolaides coarse space.  This improvement comes at the price of a somewhat larger coarse problem.
Current work in progress involves coarse approximation error and stable decomposition estimates and is left to a future article by the same authors.  
We are also planning to extend the presented two-level preconditioning strategy to nonlinear PDEs that model free-surface flows.


\section*{Acknowledgments}
This work has been supported by ANR Project Top-up (ANR-20-CE46-0005).  
The high-resolution structural data has been provided by Métropole
Nice Côte d'Azur.  We warmly thank Florient Largeron,  chief of MNCA's
SIG 3D project,  for his help in preparation of the data and for the multiple fruitful
discussions.

\bibliography{biblioabbrev}

\begin{thebibliography}{10}
\providecommand{\url}[1]{{#1}}
\providecommand{\urlprefix}{URL }
\expandafter\ifx\csname urlstyle\endcsname\relax
  \providecommand{\doi}[1]{DOI~\discretionary{}{}{}#1}\else
  \providecommand{\doi}{DOI~\discretionary{}{}{}\begingroup
  \urlstyle{rm}\Url}\fi

\bibitem{multiscalehighaspectratio}
Aarnes, J., Hou, T.Y.: Multiscale domain decomposition methods for elliptic
  problems with high aspect ratios.
\newblock Acta Math. Appl. Sin. Engl. Ser. \textbf{18}(1), 63--76 (2002)

\bibitem{diffwave}
Alonso, R., Santillana, M., Dawson, C.: On the diffusive wave approximation of
  the shallow water equations.
\newblock Eur. J. Appl. Math. \textbf{19}(5), 575--606 (2008)

\bibitem{buildings}
Andres, L.: L’apport de la donn{\'e}e topographique pour la mod{\'e}lisation
  3{D} fine et classifi{\'e}e d’un territoire.
\newblock Rev. XYZ \textbf{133}(4), 24--30 (2012)

\bibitem{taralova}
Brown, D.L., Taralova, V.: A multiscale finite element method for {N}eumann
  problems in porous microstructures.
\newblock Discrete \& Continuous Dynamical Systems-S \textbf{9}(5), 1299 (2016)

\bibitem{gdsw}
Dohrmann, C.R., Klawonn, A., Widlund, O.B.: Domain decomposition for less
  regular subdomains: Overlapping {S}chwarz in two dimensions.
\newblock SIAM J. Numer. Anal. \textbf{46}(4), 2153--2168 (2008)

\bibitem{multileveldiscon}
Dryja, M., Sarkis, M.V., Widlund, O.B.: Multilevel schwarz methods for elliptic
  problems with discontinuous coefficients in three dimensions.
\newblock Numer. Math. \textbf{72}(3), 313--348 (1996)

\bibitem{analysis}
Dryja, M., Smith, B.F., Widlund, O.B.: Schwarz analysis of iterative
  substructuring algorithms for elliptic problems in three dimensions.
\newblock SIAM J. Numer. Anal. \textbf{31}(6), 1662--1694 (1994)

\bibitem{eigenproborig}
Galvis, J., Efendiev, Y.: Domain decomposition preconditioners for multiscale
  flows in high-contrast media.
\newblock Multiscale Model. Simul. \textbf{8}(4), 1461--1483 (2010)

\bibitem{shemorig}
Gander, M.J., Loneland, A., Rahman, T.: Analysis of a new harmonically enriched
  multiscale coarse space for domain decomposition methods.
\newblock arXiv preprint arXiv:1512.05285  (2015)

\bibitem{multiscalepdes}
Graham, I.G., Lechner, P., Scheichl, R.: Domain decomposition for multiscale
  {PDE}s.
\newblock Numer. Math. \textbf{106}(4), 589--626 (2007)

\bibitem{legoll}
Le~Bris, C., Legoll, F., Lozinski, A.: An {M}s{FEM} type approach for
  perforated domains.
\newblock Multiscale Model. Simul. \textbf{12}(3), 1046--1077 (2014)

\bibitem{balancing}
Mandel, J., Brezina, M.: Balancing domain decomposition for problems with large
  jumps in coefficients.
\newblock Math. Comput. \textbf{65}(216), 1387--1401 (1996)

\bibitem{dtn}
Nataf, F., Xiang, H., Dolean, V., Spillane, N.: A coarse space construction
  based on local {D}irichlet-to-{N}eumann maps.
\newblock SIAM J. Sci. Comput. \textbf{33}(4), 1623--1642 (2011)

\bibitem{nic}
Nicolaides, R.A.: Deflation of conjugate gradients with applications to
  boundary value problems.
\newblock SIAM J. Numer. Anal. \textbf{24}(2), 355--365 (1987)

\bibitem{geneo}
Spillane, N., Dolean, V., Hauret, P., Nataf, F., Pechstein, C., Scheichl, R.:
  Abstract robust coarse spaces for systems of {PDE}s via generalized
  eigenproblems in the overlaps.
\newblock Numer. Math. \textbf{126}(4), 741--770 (2014)

\bibitem{bemfem}
Wei{\ss}er, S.: BEM-based Finite Element Approaches on Polytopal Meshes, vol.
  130.
\newblock Springer (2019)

\end{thebibliography}
\bibliographystyle{spmpsci}
\end{document}